\documentclass[10pt]{amsart}
\usepackage{latexsym, amscd, amsmath, verbatim, amssymb, graphics}
\theoremstyle{plain}
\newtheorem{theorem}{Theorem}[section]
\newtheorem{corollary}[theorem]{Corollary}
\newtheorem{proposition}[theorem]{Proposition}
\newtheorem{lemma}[theorem]{Lemma}
\newtheorem*{theorem3.1}{Theorem 3.1}
\newtheorem{claim}{Claim}
\newtheorem{fact}{Fact}
\theoremstyle{remark}
\newtheorem{remark}[theorem]{Remark}
\theoremstyle{definition}
\newtheorem{definition}[theorem]{Definition}

\numberwithin{equation}{theorem}

\newcommand{\In}{\operatorname{in}}

\newcommand{\tor}{\operatorname{Tor}}
\newcommand{\ext}{\operatorname{Ext}}

\begin{document}

\title{The pinched Veronese is Koszul}
\author{Giulio Caviglia}
\date{\today}
\begin{abstract}
In this paper we prove that the coordinate ring of the pinched
Veronese (i.e $k[X^3,X^2Y,XY^2,Y^3,X^2Z,Y^2Z,XZ^2,YZ^2,Z^3]\subset
k[X,Y,Z]$) is Koszul. The strategy of the proof is the following: we
can consider a presentation $S/I$ where $S=k[X_1,\dots,X_9]$. Using a
distinguished weight $\omega$, it's enough to show that
$S/in_{\omega}I$ is Koszul. We write $in_{\omega}I$ as $J+H$ where $J$
is generated by a Gr\"obner basis of quadrics. Finally, we present an
extension of the notion of Koszul filtration and we use it to show
that $(J+H)/J$ has a linear free resolution over $S/J.$ This implies the Koszulness of $S/I.$
\end{abstract}
\address{Mathematics Department\\
         University of Kansas \\
         405 Snow Hall\\
         1460 Jayhawk Blvd\\
         Lawrence, KS 66045-7523, USA}
\email{caviglia@math.ukans.edu}
\thanks{The author was partially supported by the ``Istituto Nazionale
di Alta Matematica Francesco Severi'', Rome.}
\date{}

\maketitle

\section{Introduction}
A standard graded $k$-algebra $R$ is said to be {\it Koszul} if its residue
field $k$ has a linear free resolution as an $R$-module. It is a well
known fact that a Koszul algebra $R$ has to be quadratic in the sense
that there exists a presentation $R\cong k[X_1,\dots,X_n]/I$ where $I$
is generated by homogeneous forms of degree two. The converse does
not hold in general; the first counterexample was found by C.Lech and
consists of $k[X_1,\dots,X_4]/I$ where $I$ is generated by five generic
quadratic forms. On the other hand the result is true for example
when $I$ is given by a complete intersection of quadrics (in this case the result can be
deduced by the knowledge of the Poincare series of $R/I$, determined
in \cite{T}) or when $I$ is generated by a
Gr\"obner basis, with respect to some term order, consisting of
elements of degree two (see \cite{A}).
Note that both the previous conditions are sufficient to imply
Koszulness but not
necessary. Several results concerning Koszul algebras can be found in \cite{F}.

Many classical varieties, like
Grassmannians, Schubert varieties, flag varieties are not only Koszul but
they are presented by a quadratic Gr\"obner basis in their natural
embedding. For example Kempf \cite{K} proved that the coordinate ring of at
most $2n$ points of $\mathbb{P}^n$ in general position  is Koszul and
later A.Conca, N.V.Trung and G.Valla \cite{CTV} showed that it also admits a
quadratic Gr\"obner basis.
Another important aspect is the relation between the Koszulness
(and more generally the study of the free resolution of a residue field)
and the structure of the non commutative algebra $\ext^{*}_R(k,k)$,
i.e the Yoneda-Hopf algebra of $k.$

The study of the Koszulness for semi-group rings (i.e. toric varieties)
and the connection with the related combinatorial objects have given
the motivation for the development of interesting results (see for
example \cite{BGT}, \cite{HHR} and \cite{S2}).

Given $R$ any commutative graded $k$-algebra the Veronese subalgebra
$R^{(d)}$ has a quadratic Gr\"obner basis for $d \gg 0$ (\cite{ERT})
and so in particular it is Koszul.
This result was later generalized to a larger
class of algebras not necessarily generated in degree one, see
\cite{BGT}. Sturmfels in \cite{S2} has shown that the $k$-algebra
called of {\it Veronese type} i.e the subring of $S=k[x_1,\dots,x_n]$
generated by the monomials $\{x_1^{i_1}\cdots x_n^{i_n} | i_1+\dots
+i_n=d, 0\leq i_1\leq s_1,\dots, 0\leq i_n\leq s_n \}$ has a Gr\"obner
basis, in a certain ordering, which is not only quadratic but also
square-free. Note that when $s_1=\dots =s_n=d$ we get $S^{(d)}.$
Further generalizations have been done by S.Blum in \cite{B}.

An important question, regarding the Koszulness of toric variety,
which is still open is the following: ``Is it true that any quadratic
toric varieties with an isolated singularity is Koszul?''
The pinched Veronese, i.e. the $k$-algebra, where $k$ is a field,
defined as $R=k[X^3,X^2Y,XY^2,Y^3,X^2Z,Y^2Z,XZ^2,YZ^2,Z^3],$ has
been for a long time the first and the most simple case of the previous
question where the answer was unknown. The problem about the
Koszulness of the pinched Veronese was raised by B.Sturmfels in the
1993 in a conversation with Irena Peeva, and after that has been circulating
as a concrete example to test the efficiency
of the new theorems and techniques concerning Koszul algebras.

In Section 3 of this paper we are going to prove that the pinched
Veronese is Koszul, i.e.
\begin{theorem3.1} The algebra
$R=k[X^3,X^2Y,XY^2,Y^3,X^2Z,Y^2Z,XZ^2,YZ^2,Z^3],$ where $k$ is a field, is Koszul.
\end{theorem3.1}

\section {Initial ideal with respect to a weight and Koszul
filtrations for modules}

It is well known that a graded algebra $k[X_1\ldots,X_n]/I$ with $I\subseteq (X_1,\dots ,X_n)^2$ is Koszul if, for a certain term order, it can be generated by a Gr\"obner basis of
quadratic forms. If $I$ is monomial the converse holds \cite {BHV} , consequently, if for a certain term order $\In(I)$ defines a Koszul algebra, so does $I$.
It is known that the same result is still true
considering, instead of a term order, a weight function $w$ given by a
vector $(w_1, \ldots, w_n)$ of positive integers and  replacing
$\In(I)$ by the initial ideal $\In_{w}(I)$ (not necessary monomial) of $I$ with
respect to $w.$ The proof of this fact is based on the idea of flat
families and we refer for notations and generalities concerning them to \cite{E} 15.8.

Given a weight function $w=(w_1,\ldots,w_n)$ from $\mathbb{Z}^n$ to
$\mathbb{Z}$ we can think about it as a function defined on monomials
of $R$; moreover given $f\in R$ we use $\In_w(f)$ for the sum of all
the terms of $f$ that are maximal with respect to $w.$ Given an
ideal $I$ we write $\In_w(I)$ for  the ideal generated by
$\In_w(f)$ for all $f\in I.$
Let $S=R[T]$ be a polynomial ring in one variable over $R;$ for any
$f\in R$ we define $\tilde{f}$ in the following way: we can write
$f=\sum u_im_i$ where $m_i$ are distinct monomials and $0\not= u_i \in k.$ Let
$a=\max w(m_i),$ and set

$$\tilde{f}=T^af(T^{-w_1}X_1,\dots,T^{-w_n}X_n).$$

Note that $\tilde{f}$ can be written as $\In_{\omega}(f)+gT$ where $g$ belongs
to $S.$ For any ideal $I$ of $R$ define $\tilde{I}$ to be the ideal
of $S$ generated by the elements $\tilde{f}$ for all $f\in I.$
Setting $\deg X_i=(1,w_i)$ and $\deg T=(0,1)$ the algebra $S$ is
bigraded and in particular if $I\subset R$ is an homogeneous ideal then
$\tilde{I}$ is bihomogeneous. From the definition it follows that
$S/((T)+\tilde{I})\cong R/\In(I)$ and $S/((T-1)+\tilde{I})\cong
R/I.$ Note that $T$ is a non-zerodivisor on $S/\tilde{I}:$ let
$Tf\in \tilde{I}$ for some $f\in S.$ Without loss of generality we
can assume $f$ bihomogeneous and moreover specializing at $T=1$ we have
$h=f(X_1,\ldots,X_n,1)\in I,$ but $f$ is bihomogeneous therefore it
holds that $f=\tilde{h}\in \tilde{I}.$ From this fact it follows that also
$T-1$ is a non-zerodivisor for $S/\tilde{I}$ since it is sum of two
non-zerodivisors of different degrees.

The following statement is folklore, but since we couldn't find a proof in the literature we give one here for completeness.

\begin{lemma} \label{lemma1} Let $R=k[X_1,\ldots,X_n]$. Consider a weight given by a vector of positive integers $w=(w_1,\ldots,w_n)$ and homogeneous ideals $I,J,H$
such that $I\subseteq J$ and $I\subseteq H.$ Then
$$\dim_k \tor_i^{R/I}(R/J,R/H)_j\leq \dim_k \tor
_i^{R/\In_w I }(R/\In_w J, R/\In_w H)_j.$$
\end{lemma}

\begin{proof}  Consider the ideals $\tilde {I},$
$\tilde {J},$ $\tilde {H}$ of $S=R[T]$ defined as above. Let
$M_{i}=\tor_i^{S/\tilde{I}}(S/\tilde{J},S/\tilde{H})$; note that $M_i$ is
bigraded and we can make it a $\mathbb{Z}$-graded module setting
$(M_i)_j=\bigoplus _{h\in \mathbb{Z}}(M_i)_{(j,h)}.$ Since $(M_i)_j$ is a
finitely generated module over $k[T]$, the structure theorem for
modules over a PID
applies and we obtain the isomorphism $(M_i)_j\cong k[T]^{a_{ij}} \bigoplus A_{ij}$ where $A_{ij}$ is the torsion submodule. Moreover since $A_{ij}$ has to be homogeneous the
structure theorem gives $A_{ij}\cong \bigoplus_{h=1} ^{b_{ij}}k[T]/(T^{d_h}).$
Set $l_1=T$ and $l_2=T-1$ and consider the following exact sequence
\begin {eqnarray} \label{aa}
\begin{CD} 0 @>>> S/\tilde{H} @>\cdot l_r>>
S/\tilde{H} @>>> S/((l_r)+\tilde{H}) @>>> 0,
\end{CD}
\end {eqnarray}
for $r=1,2.$
All the modules appearing in \eqref{aa} are over $S/\tilde{I}$ and the
multiplication by $l_i$ is a zero degree map with respect to the
$\mathbb {Z}-$grading.
Tensoring with $S/\tilde{J}$ and passing to the long exact sequence of
homologies we have
\begin{equation}\label{bo} 0\rightarrow M_i/l_r M_i \rightarrow
\tor_i^{S/\tilde{I}}(S/\tilde{J},S/((l_r)+\tilde{H})) \rightarrow
\ker (M_{i-1} \xrightarrow{\cdot l_r} M_{i-1})\rightarrow 0.
\end{equation}

Since $l_i$ is a regular element for $S/\tilde{I}$ and $S/\tilde{J}$,
the middle term is isomorphic to
$\tor_i^{S/((l_r)+\tilde{I})}(S/((l_r)+\tilde{J}),S/((l_r)+\tilde{H}))$
(see \cite{M} Lemma 2 page 140)
which is $\tor_i^{R/\In_w{I}}(R/\In_w{J},R/\In_w{H}))$ when $r=1$ and is $\tor_i^{R/I}(R/J,R/H))$ for $r=2.$

Therefore taking the graded component of degree $j$ in \eqref{bo}  we obtain:
\begin{eqnarray} \label{do}
\dim_k
\tor_i^{R/\In_w{I}}(R/\In_w{J},R/\In_w{H})_j &=& a_{ij}+b_{ij}+b_{(i-1)j}, \\ \label{el}
\dim_k
\tor_i^{R/I}(R/J,R/H))_j &=& a_{ij}.
\end{eqnarray}

The Lemma follows by comparing \eqref{do} and \eqref{el}.
\end{proof}

\begin{corollary}\label{corollary1} Let $w$ be a weight function and $I$ be a homogeneous ideal of
$R=k[X_1,\ldots,X_n]$ such that $R/\In_w(I)$ is Koszul. Then $R/I$ is
Koszul.
\end{corollary}

\begin{proof} Since $R/\In_w I $ is Koszul, $\dim_k
\tor_i^{R/\In_w I }(k,k)_j=0$ for any $i\not= j.$ Applying Lemma \ref{lemma1} with $J=H=(X_1,\ldots, X_n),$ we see that  $\tor_i^{R/I}(k,k)_j=0$ for any $i\not= j.$
\end{proof}

We recall the definition of a \emph{Koszul filtration} introduced
by A.Conca, N.V.Trung and G.Valla in \cite{CTV}, see also \cite{HHR} for related results.
\begin{definition} \label{koszul filtration} Let $R$ be a standard graded $k$-algebra. A family
$\mathbf{F}$ of ideals of $R$ is said to be a \emph{Koszul
filtration} of $R$ if:\\
1) Every ideal $I\in \mathbf{F}$ is generated by linear forms.\\
2) The ideal $(0)$ and the maximal homogeneous ideal $\mathcal {M}$ of
$R$ belong to $\mathbf{F}$.\\
3) For every $I\in \mathbf{F}$ different from $(0)$, there exists $J\in
\mathbf{F}$ such that $J\subset I$, $I/J$ is cyclic and $J:I\in
\mathbf{F}.$
\end{definition}

In \cite{CTV} it is proved that all the ideals belonging to such a
filtration have a linear free resolution over $R$ and in particular, since the
homogeneous maximal ideal is in $\mathbf{F}$, $R$ will be a Koszul algebra.
For the purposes of this paper, we need an extension of
this definition to the case of graded modules. In particular we want
$\mathbf{F}$ to be a collection of finitely generated graded modules.

\begin{definition}\label{generalized Koszul filtration} Let $R$ be a
standard graded $k$-algebra. A family $\mathbf{F}$ of finitely
generated graded $R$-modules is said to be a \emph{Koszul filtration}
for modules if the following three properties hold:\\
1) Every nonzero module $M \in \mathbf{F}$ is generated by its nonzero component of lowest degree, say $s_M.$\\
2) The zero module belongs to $\mathbf{F}.$\\
3) For every $M\in \mathbf{F}$ different from the zero module there
exists   $N\subsetneq M,$ $N\in \mathbf{F}$ with $N=0$ or $s_M=s_N$,
such that either $M/N$ has a linear free resolution (i.e
$\tor_i^R(M/N,k)_j=0$ for all $j\not =i+s_M$) or the module of first
syzygies $\Omega_1^R(M/N)$ of $M/N$ is generated in degree $s_M+1$ and $\Omega_1^R(M/N)(1) \in
\mathbf{F}.$
\end{definition}

The next proposition shows that all the elements in $\mathbf{F}$ have a
linear free resolution over R. In particular, because of this fact,
the problem of proving that some module $M$ has a linear free resolution
over $R$ can be approached by trying to construct a \emph{Koszul
filtration} containing $M.$ This way of proceeding is
very useful in the next section where the \emph{Koszul filtration}
for modules plays an important role in proving that the pinched Veronese
is Koszul.

\begin{proposition}\label{lemma2} Let $R$ be a standard graded $k$-algebra and
$\mathbf{F}$ a Koszul filtration as defined in \eqref{generalized Koszul filtration}.
Then every $M\in \mathbf{F}$ has a linear free resolution over $R.$
\end{proposition}
The proof of this result is essentially the same as the one for the case
of a Koszul filtration (see \cite{CTV} {Prop. 1.2}).
\begin{proof} We need to show that for every $M\in \mathbf{F}$ we have
$\tor_i^R(M,k)_j=0$ for all $j\not= i+s_M.$ We argue by induction on the index $i.$
If $i=0$ the assertion is clearly true; fix an integer $i>0.$ If $M$
is the zero module it has obviously a linear resolution, therefore we
can assume that $M$ has a positive minimum number of generators $\mu
(M).$ Inducting on $\mu(M)$ we can assume that $\tor_i^R(N,k)_j=0$ whenever $j\not= i+s_N,$  $N\in \mathbf{F}$ and $\mu(N)<\mu(M)$.
The third property in Definition \ref{generalized Koszul filtration} implies that $M$ has a submodule $N\subsetneq M$ with
$s_M=s_N$ and in particular with $\mu(N)<\mu(M).$\\
The short exact sequence
$$0 \longrightarrow N \longrightarrow M \longrightarrow M/N
\longrightarrow 0 $$
gives the exact sequence
\begin{equation} \label{ex seq tor}
\tor_i^R(N,k)_j \longrightarrow \tor_i^R(M,k)_j \longrightarrow
\tor_i^R(M/N,k)_j .
\end{equation}
>From the third property in Definition \ref{generalized Koszul filtration} either $M/N$ has a linear resolution, so in particular $\tor_i^R(M/N,k)_j=0$ for $j\not=i+s_M,$ or $\Omega_1^R(M/N)(1)\in \mathbf{F}$ and is generated in degree $s_M.$
The last term in \eqref{ex seq tor} is isomorphic to
$\tor_{i-1}^R(\Omega_1^R(M/N)(1),k)_{j-1}.$
Since the inductive hypothesis on the index of the $\tor$ applies we
deduce that $\tor_{i}^R(M/N,k)_j=\tor_{i-1}^R(\Omega_1^R(M/N)(1),k)_{j-1}= 0$ when $j\not=i+s_M$. On the other hand the induction on the minimum number of generators yields
$\tor_i^R(N,k)_j=0$ for $j\not=s_M+i$ and therefore the middle
term in \eqref{ex seq tor} vanishes when $ j\not=s_M+i.$
\end{proof}

\begin{remark} We consider, in the third property in Definition \ref{generalized
Koszul filtration}, the fact of having a linear free resolution over $R$ as a possible condition  for an element $M$ in $\mathbf{F}$. This
is not really essential: indeed if we
already know that $M$ has a linear free resolution over $R$ we could add to
$\mathbf {F}$ all the
modules of syzygies of $M$ filtering them, trivially, with $0$. We decide, for the sake of
convenience, to try to keep the family $\mathbf{F}$ as small as possible.
On the other hand, if we only leave the second part of condition number
3), it reasonable to ask,
given a module $M$ with a linear free resolution over $R$, if there
always exists a finite family $\mathbf{F}$ containing $M$.
\end{remark}

\begin{remark} It's easy to see that the Koszul filtration is included
in our definition of Koszul filtration for modules.
In fact if $J\subset I$ are ideals generated by linear forms (as in Definition
\ref{koszul filtration}) with $I/J$ cyclic, then $J:I \cong
\Omega_1^R(I/J)(1).$
\end{remark}
\begin{remark}
Our definition of Koszul filtration covers also the definition of
\emph{module with linear quotients} recently introduced by Conca and
Herzog in \cite{CH} in order to study the linearity of the free
resolution of certain modules over a polynomial ring.
\end{remark}

\section{The pinched Veronese}
In this section we prove that the coordinate ring $R$ of the pinched
Veronese variety is a Koszul algebra. The proof is structured in three
different steps.
First of all we can consider a presentation for $R$ and write it as
$R=S/I$ where $I$ is a homogeneous ideal generated by quadrics and
$S$ is a polynomial ring.

The first step consists in taking the initial ideal of $I$ with
respect to a carefully chosen weight $\omega.$ By Corollary
\ref{corollary1}, it's then sufficient to show that $S/\In_{\omega}(I)$ is
Koszul. The use of $\omega$ is important because it allows us to study
instead of a binomial ideal, an ideal generated by
several quadratic monomials and only five quadratic binomials. For this
purpose the choice of $\omega$ needs to be done very carefully:
taking, for example, the trivial weight $\omega=(1,\dots,1)$ we get
$\In_{\omega}(I)=I$ and we do not make any simplification. On the
other hand a generic weight will play the role of a term order,
bringing in the initial ideal some minimal generator of degree higher than two,
and so the ring defined by the initial ideal with respect to a generic
weight cannot be Koszul.

The second reduction consists in writing $\In_{\omega}(I)$ as the sum
two ideals: $U$ generated by the monomial part of  $\In_{\omega}(I)$ plus a
distinguished binomial of $\In_{\omega}(I)$ and the ideal $(Q_1,\dots
Q_4)$ given  by
the remaining four binomials of  $\In_{\omega}(I).$ The ideal $U$
is generated by a Gr\"obner basis of quadrics, so $S/U$ is Koszul.
We need the following fact, which is part of Lemma 6.6 of \cite{CHTV}:

\begin{fact}\label{fact} Let $T$ be a Koszul algebra let $Q\subset T$ be a quadratic ideal
with a linear free resolution over $T.$ Then $T/Q$ is Koszul.
\end{fact}
Using this fact it is enough to show that the class $(q_1,\dots,q_4)$ of
$(Q_1,\dots,Q_4)$ in $S/U$ has a linear free resolution over $S/U.$
Note that among all the possible five binomials, the distinguished one we
pick is the only one giving at the same time the Koszulness of $S/U$
and the linearity of the ideal given by the other four. It's maybe
possible to show that the whole binomial part has a linear resolution
over $S$ modulo the monomial one, but for this purpose the amount of
calculations required seems much higher.

The last part of the proof consists in showing the linearity of the
free resolution of $(q_1,\dots,q_4)$ over $S/U$ via the construction
of a Koszul filtration containing $(q_1,\dots,q_4)$.

\begin{theorem} The algebra
$R=k[X^3,X^2Y,XY^2,Y^3,X^2Z,Y^2Z,XZ^2,YZ^2,Z^3]$, where $k$ is field, is Koszul.
\end{theorem}
\begin{proof}
Since $R$ contains all monomials in $X,Y,Z$ of degree 9 its Hilbert
function is $H_R(0)=1$, $H_R(1)=9$ and $H_R(n)=\binom{3n+2}{2}$ for
$n\geq 2.$ The Hilbert polynomial of $R$ is given by $\binom{3n+2}{2}$
and the Krull dimension of $R$ is $3$.
One computes that the Hilbert series of $R$ is given by:
$$H_R(Z)=\frac{Z^4-3Z^3+4Z^2+6Z+1}{(1-Z)^3}.$$

Consider a presentation $S/\ker \phi \cong R$ where $S=k[X_1,\dots,X_9]$ and
and $\phi$ is the homomorphism from $S$ to $R$ defined by sending $X_i$
to the $i^{th}$ monomial of $R$ in
$(X^3,X^2Y,XY^2,Y^3,X^2Z,Y^2Z,XZ^2,YZ^2,Z^3).$ Let $I$ be the ideal
defined as
\begin{eqnarray*}
&I=&(X_8^2-X_6X_9,X_6X_8-X_4X_9,X_5X_8-X_2X_9,X_7^2-X_5X_9,\\
&&X_6X_7-X_3X_9,X_5X_7-X_1X_9,X_4X_7-X_3X_8,X_3X_7-X_2X_8,\nonumber\\
&&X_2X_7-X_1X_8,X_6^2-X_4X_8,X_5X_6-X_2X_8,X_5^2-X_1X_7,X_4X_5-X_2X_6,\nonumber\\
&&X_3X_5-X_1X_6,X_3^2-X_2X_4,X_2X_3-X_1X_4,X_2^2-X_1X_3).\nonumber
\end{eqnarray*}
It is immediate to see that $I\subseteq \ker \phi.$ On the other hand
also the opposite inclusion holds, in fact it is sufficient to check that $R$ and
$S/I$ have the same Hilbert function. We will prove this below.

Consider the weight function $\omega$ from $\mathbb{Z}^9$ to $\mathbb{Z}$ given by $(3,3,1,3,3,3,2,3,3)$ and take its natural extension to the monomials of $S.$
Let $J$ be the ideal generated by the initial forms with respect to $\omega$ of the generators
of $I$ given previously. We have
\begin{eqnarray*}J&=&(X_8^2-X_6X_9,X_6X_8-X_4X_9,X_5X_8-X_2X_9,X_5X_9,X_6X_7,\\
&& X_1X_9,X_4X_7,X_2X_8,X_1X_8,X_6^2-X_4X_8,X_5X_6,X_5^2,\\
&&X_4X_5-X_2X_6,X_1X_6,X_2X_4,X_1X_4,X_2^2).
\end{eqnarray*}
We claim that $J= \In_{\omega}I$: one inclusion is clear and to prove
the other is enough to show, as stated previously, that $R/J$ and $R/I$ have
the same Hilbert function. Consider the degrevlex order $\sigma$ on the
monomials of $S$. Note first that $X_2X_9^2$ and $X_2X_6X_9$ belong to
$J$ since $X_2X_9^2=(X_5X_9)X_8-(X_5X_8-X_2X_9)X_9$ and
$X_2X_6X_9=(X_2X_8)X_8-(X_8^2-X_6X_9)X_2,$ therefore the
following ideal
\begin{eqnarray*}H&=&(X_5X_9,X_1X_9,X_8^2,X_6X_8,X_5X_8,
X_2X_8,X_1X_8,X_6X_7,X_4X_7,X_6^2,\\
&&X_5X_6,X_1X_6,X_5^2,X_4X_5,X_2X_4,X_1X_4,X_2^2,X_2X_9^2,X_2X_6X_9)\end{eqnarray*}
is contained in $\In_{\sigma}J.$ The Hilbert series of $S/H$ is easy
to compute and it is
$H_{S/H}(Z)=\frac{Z^4-3Z^3+4Z^2+6Z+1}{(1-Z)^3}$. Coefficient-wise we
have:
\begin{eqnarray*}\label{hil} H_{S/\ker \phi}(Z)\leq H_{S/I}(Z)= H_{S/\In_{\omega}I}(Z) \leq H_{S/J}(Z)=
H_{S/\In_{\sigma}J}(Z)\leq H_{S/H}(Z).
\end{eqnarray*}
The first and the last term agree, thus all the previous
inequalities are in fact equalities, and in particular it
follows that $\ker \phi=I,$ $\In_{\omega}I=J$ and $\In_{\sigma}J=H.$

Applying Corollary \ref{corollary1} to $S/I$ and $\omega$, in order to
finish the proof of the theorem it's
enough to show the following
\begin{claim} \label{claimA}
The $k$-algebra $S/J$ is Koszul.
\end{claim}
\emph{Proof of the claim }
We can write $J$ as a sum of two ideals: one generated by all the
quadratic monomials of $J$ together with the quadratic binomial
$X_6X_8-X_4X_9,$ namely
\begin{eqnarray*}
U&=&(X_5X_9,X_1X_9,X_2X_8,X_1X_8,X_6X_7,X_4X_7,X_5X_6,X_1X_6,X_5^2,X_2X_4,\\
&&X_1X_4,X_2^2,X_6X_8-X_4X_9),
\end{eqnarray*}
and the other one generated by the remaining binomials $Q_1=X_6^2-X_4X_8$, $Q_2=X_4X_5-X_2X_6$, $Q_3=X_8^2-X_6X_9$ and
$Q_4=X_5X_8-X_2X_9.$
Note first that $U$ is generated by a Gr\"obner basis of quadrics with
respect to the degrevlex order $\sigma,$ in fact all the S-pairs we need to
check are:
$$(X_6X_8-X_4X_9)X_2-(X_2X_8)X_6=-(X_2X_4)X_9,$$
$$(X_6X_8-X_4X_9)X_1-(X_1X_8)X_6=-(X_1X_9)X_4,$$
\begin{equation}\label{grobner}
(X_6X_8-X_4X_9)X_7-(X_6X_7)X_8=-(X_4X_7)X_9,
\end{equation}
$$(X_6X_8-X_4X_9)X_5-(X_5X_6)X_8=-(X_5X_9)X_4,$$
$$(X_6X_8-X_4X_9)X_1-(X_1X_8)X_6=-(X_1X_9)X_4.$$
One can observe that for any ideal
$L=(X_{i_1},\dots,X_{i_r})$ generated by variables, the ideal $U+L$ is again generated by a
Gr\"obner basis of at most quadrics. Indeed there are no S-pairs
to check other than the ones in \eqref{grobner}. Moreover if $L$ is chosen in
a such a way that  $X_6X_8-X_4X_9 \in L$ or $X_6X_8\not
\in L,$ we obtain that  $\In_{\sigma}(U)
+L=\In_{\sigma}(U+L),$ in fact the only case in which this very last
equality doesn't hold is when $X_4X_9$ appears in the sum without being
in $L.$
By Theorem 2.2 of \cite{BHV} if $\In_{\sigma}(U)
+L=\In_{\sigma}(U+L)$ then not only $S/U$ is Koszul but also
the ideal $(L+U)/U$ has a linear free resolution over
$S/U.$

Now set $S/U=T$. In the following we will denote by $x_i$ the class of $X_i$ and by
$q_j$ the class of $Q_j$ in $T.$
Since $T$ is Koszul we can use Fact \ref{fact} to conclude the proof
of Claim \ref{claimA} by showing that $(q_1,\dots,q_4)$ has a linear free resolution over $T.$
We prove this by constructing  a
\emph{Koszul filtration} $\mathbf{F}$ over $T$ containing
$(q_1,\dots,q_4)(1)$ because this implies, by Proposition \ref{lemma2}, that
$(q_1,\dots,q_4)(1)$ has a linear free resolution over $T$
and so $(q_1,\dots,q_4)$ does.

It will be useful to include in $\mathbf{F}$ a set of ideals $\mathbf{G}$
for which we already know they have a linear free resolution over $T.$
Setting
\begin{eqnarray*}
\mathbf{G}&=&\{\mathrm{Ideals}\;(x_{i_1},\dots,x_{i_r})\mathrm{\; of\;}T \mathrm{\; such\;
that\; } X_6X_8-X_4X_9 \in (X_{i_1},\dots,X_{i_r})\\
&&\mathrm{\;\;or\;} X_6X_8\not
\in (X_{i_1},\dots,X_{i_r})\},
\end{eqnarray*}
from what we have seen above any ideal in $\mathbf G$ has a linear
resolution over $T.$

We define $\mathbf{F}$ to be
\begin{align*} \mathbf{F}&=\mathbf{G}\cup
\{(q_1,q_2)(1),(q_1,\dots,q_4)(1), M_1^{1,\dots,8},
 M_1^{1,2,4,6,7,8},\\
&M_2^{1,\dots,8},M_2^{1,2,4,5,7,8},M_3^{1,\dots,10},M_3^{1,2,3,5,6,7,8,9,10}
 \} \cup \{0\}.
\end{align*}

The modules $M_1^{1,\dots,8},$
 $M_1^{1,2,4,6,7,8},$
$M_2^{1,\dots,8},$$M_2^{1,2,4,5,7,8},$$M_3^{1,\dots,10},$$M_3^{1,2,3,5,6,7,8,9,10}$ are constructed as follows. We consider $T$-homomorphisms defined by matrices:
\begin{align*}
M_1&=\begin{pmatrix} x_7 & 0 & x_5 & x_1 & x_2 & 0 & 0 & 0 \\
                        0   &  x_7    & x_8    & 0   & x_6 & x_5 & x_2
                        & x_1
\end{pmatrix} \\
M_2&=\begin{pmatrix} x_7 &  0 & -x_8 & 0 & x_2 & -x_5 & x_1 & 0 \\
                        0   &  x_7    & x_6   & x_5 & 0 & x_2 & 0 &
                        x_1
\end{pmatrix}\\
M_3&=\begin{pmatrix} x_6 &  0 & x_5 & -x_8 & x_4 & 0& x_2 & 0 & x_1 & 0 \\
                        0   &  x_6    & 0   & x_5 & 0 & x_4 & 0 & x_2
                        & 0 & x_1
\end{pmatrix}
\end{align*}

$$ M_1: T(-1)^8 \rightarrow T\,^2,
\quad M_2: T(-1)^8 \rightarrow T\,^2,
\quad M_3: T(-1)^{10} \rightarrow T\,^2.$$

We use now an upper index notation on the matrices to indicate the module generated by
the images of the elements of the standard basis corresponding to those
indeces: for instance  $M_1^{1,4}$ is the module generated by the
images under $M_1$ of $(1,0,\dots,0)$ and $(0,0,0,1,0,\dots,0)$.

We prove that $\mathbf{F}$ is a \emph{Koszul filtration}
for $T.$
For what concerns the elements in $\mathbf{G}$ there is nothing to
check since $0\in \mathbf{F}$ and they have a linear free
resolution over $T$. For all the other modules $M \in \mathbf{F}$ we
exhibit a submodule $N\in \mathbf{F},\, N\subsetneq M$, such that
$M/N$ has a linear free resolution or $\Omega_1(M/N)(1)$ belongs to
$\mathbf{F}.$ We have the following isomorphisms
\begin{align}
\Omega_1((q_1,q_2)(1))(1)&\cong M_1^{1,\dots,8} \in \mathbf{F} \label{a}\\
\Omega_1((q_1,\dots,q_4)/(q_1,q_2)(1))(1)&\cong
M_3^{1,\dots,10}\in\mathbf{F} \label{b}\\
\Omega_1(M_1^{1,\dots,8}/M_1^{1,2,4,6,7,8})(1)&\cong M_2^{1,\dots,8}\in
\mathbf{F}\label{c}\\
M_1^{1,2,4,6,7,8}/M_1^{1,4}&\cong(x_7,x_5,x_2,x_1) \in
\mathbf{G}\subseteq \mathbf {F}  \label{d}\\
\Omega_1(M_2^{1,\dots,8}/M_2^{1,2,4,5,7,8})(1)&\cong M_1^{1,\dots,8}\in
\mathbf{F}\label{e}\\
M_2^{1,2,4,5,7,8}/M_2^{1,5,7}&\cong (x_7,x_5,x_1) \in
\mathbf{G}\subseteq \mathbf {F} \label{f}\\
\Omega_1(M_3^{1,\dots,10}/M_3^{1,2,3,5,6,7,8,9,10})(1)&\cong (x_6,x_5,x_4,x_2,x_1)\in
\mathbf{G}\subseteq \mathbf {F}\label{g}\\
M_3^{1,2,3,5,6,7,8,9,10}/M_3^{1,3,5,7,9}&\cong
(x_6,x_4,x_2,x_1) \in
\mathbf{G}\subseteq \mathbf {F} \label{h}
\end{align}
where \eqref{a},\,\eqref{b},\,\eqref{c} and \eqref{e} have been checked with the help of the computer algebra system MACAULAY2 \cite{M2} over
the field of rational numbers. In particular by flat extension these
isomorphisms work over any field of characteristic zero. On the other
hand we performed by hand exactly the same Gr\"obner basis
based computation, suggested by the calculations over ${\mathbb
Q}$. Since integer coefficients different from $1$ or $-1$ never
appear, those calculations are enough to prove the previous
isomorphisms also  over any field of positive characteristic.

In \eqref{d} and \eqref{f} the modules
$M_1^{1,4}$ and $M_2^{1,5,7}$ are clearly isomorphic to $(x_7,x_4)\in
\mathbf{G}\subseteq \mathbf{F}$ and to $(x_7,x_2,x_1)\in \mathbf{G}\subseteq \mathbf{F}$ respectively. Similarly in \eqref{h} the module $M_3^{1,3,5,7,9}$ is isomorphic to $(x_6,x_5,x_4,x_2,x_1)$
which belongs to $\mathbf{G}\subseteq \mathbf{F}$.
This show that $\mathbf {F}$ is a \emph{Koszul
filtration} and, as we said before, by Proposition \ref{lemma2} the ideal  $(q_1\dots,q_4)(1)$ has a linear free resolution over $T.$ Thus the claim is proved and so is the theorem.
\end{proof}

\section*{acknowledgements}
I would like to thank Professor Craig Huneke for many helpful comments
and several discussions which have made this paper possible.

\end{document}